\documentclass{elsart}
\usepackage{natbib}
\usepackage{amsthm, amsmath, amsfonts}

\begin{document}

\begin{frontmatter}

\title{Curve implicitization in the bivariate tensor-product Bernstein basis}

\author{Ana Marco \thanksref{EM1}}
\author{, Jos\'e-Javier Mart{\'\i}nez \thanksref{EM2}}

\address{Departamento de Matem\'aticas, Universidad
de Alcal\'a,}

\address{Campus Universitario, 28871-Alcal\'a de Henares (Madrid), Spain}

\thanks[EM1]{E-mail: ana.marco@uah.es}

\thanks[EM2]{Corresponding author. E-mail: jjavier.martinez@uah.es}

\begin{abstract}
The approach to curve implicitization through Sylvester and
B\'ezout resultant matrices and bivariate interpolation in the
usual power basis is extended to the case of Bernstein-Bezoutian
matrices constructed when the polynomials are given in the
Bernstein basis. The coefficients of the implicit equation are
also computed in the bivariate tensor-product Bernstein basis, and
their computation involves the bidiagonal factorization of the
inverses of certain totally positive matrices.


\end{abstract}

\end{frontmatter}

\section{Introduction}

When studying rational plane algebraic curves, there are two
standard ways of representation, the implicit equations and the
parametric equations. The intersection of two curves is more
easily computed when we have the implicit equation of one curve
and the parametric equations of the other, and hence it is very
important to be able to change from one representation to another.

We will concentrate on the implicitization problem, that is to
say, on finding an implicit representation starting from a given
rational parametrization of the curve.

In [12] we have presented an approach to the implicitization
problem based on interpolation using the usual power basis for the
corresponding space of bivariate polynomials. However, very recent
work [2], related to polynomials expressed in the Bernstein basis
has showed the importance of evaluating resultants from Bernstein
basis resultant matrices directly, avoiding a basis
transformation. In this sense, in [2] it is indicated that {\it
for numerical computations involving polynomials in Bernstein form
it is essential to consider algorithms which express all
intermediate results using this form only}.

Although those papers study univariate polynomials, it must be
observed that the construction of the resultant matrices can be
extended to the case in which the entries of the resultant matrix
are polynomials. It must also be taken into account that the
Bernstein basis has also important advantages in the context of
tracing implicit algebraic curves [13].

So our aim is to use bivariate interpolation for obtaining in the
bivariate tensor-product Bernstein basis the implicit equation of
a plane algebraic curve given by its parametric equations in
Bernstein form (which is the usual situation in the case of
B\'ezier curves). Although we present all the details with an
example in exact rational arithmetic, it must be taken into
account that the process can also be carried out in (high) finite
precision arithmetic. In that situation some important results of
{\sl numerical linear algebra} we use will have a major
importance. More precisely the {\it total positivity} of certain
matrices will be an important issue, as it happens in several
instances of computer aided geometric design (see, for example,
the recent work [18] and references therein).

\medskip

The rest of the paper is organized as follows. In Section $2$
several basic results will be presented. In Section $3$ we
introduce the interpolation algorithm for computing the implicit
equation as a factor of the determinant of the resultant matrix,
while in Section $4$ we consider some results related to total
positivity which will be relevant for the solving of the linear
system associated with the interpolation problem. Finally, in
Section $5$ we briefly examine the computational complexity of the
whole algorithm.

\section{Preliminaries}

Let $P(t) = (x(t),y(t))$ be a proper parametrization of a rational
plane algebraic curve $C$, where $x(t) = {{u_1(t)}\over{v_1(t)}}$
and $y(t) = {{u_2(t)}\over{v_2(t)}}$ and $gcd(u_1,v_1) =
gcd(u_2,v_2) = 1$. A parametrization $P(t) = (x(t),y(t))$ of a
curve $C$ is said to be {\sl proper} if every point on $C$ except
a finite number of exceptional points is generated by exactly one
value of the parameter $t$. It is well known that every rational
curve has a proper parametrization, so we can assume that the
parametrization is proper. Several recent results on the
properness of curve parametrizations can be seen in [17].

In connection with the implicitization problem, the following
theorem [17] holds:

\bigskip

{\bf Theorem 1}. Let
$P=(x(t)={{u_1(t)}\over{v_1(t)}},y(t)={{u_2(t)}\over{v_2(t)}})$ be
a proper rational parametri-zation of an irreducible curve $C$,
with $gcd(u_1,v_1) = gcd(u_2,v_2) = 1$. Then the polynomial
defining $C$ is  $Res_t(u_1(t) - xv_1(t), u_2(t) - yv_2(t))$ (the
resultant with respect to t of the polynomials $u_1(t) - xv_1(t)$
and $u_2(t) - yv_2(t)$).

\bigskip

Our aim is to compute the implicit equation $F(x,y) = 0$ of the
curve $C$ by means of polynomial interpolation, which taking into
account Theorem $1$ is equivalent to compute $Res_t(u_1(t) -
xv_1(t), u_2(t) - yv_2(t))$.

First of all, we remark that the concept of interpolation space
will be essential. The following result, also in [17], shows which
is in our case the most suitable interpolation space:

\bigskip

{\bf Theorem 2}. Let $P = (x(t) = {{u_1(t)}\over{v_1(t)}},y(t) =
{{u_2(t)}\over{v_2(t)}})$ be a proper rational parametri-zation of
the irreducible curve $C$ defined by $F(x,y)$, and let
$gcd(u_1,v_1) = gcd(u_2,v_2) = 1$. Then
$deg_y(F)=max\{deg_t(u_1),deg_t(v_1)\}$ and $deg_x(F)=$ $=
max\{deg_t(u_2),deg_t(v_2)\}$.

\bigskip

Theorem $2$ tells us that the polynomial $F(x,y)$ defining the
implicit equation of the curve $C$ belongs to the polynomial space
$\Pi_{n,m}(x,y)$, where $n = max \{deg_t(u_2),deg_t(v_2)\}$ and $m
= max \{deg_t(u_1),deg_t(v_1)\}$. The dimension of
$\Pi_{n,m}(x,y)$ is $(n+1)(m+1)$, and a basis is given by $\{x^i
y^j | i = 0, \cdots, n; j = 0, \cdots, m\}$. Moreover
$deg_x(F(x,y)) = n$ and $deg_y(F(x,y)) = m$, and therefore there
is no interpolation space $\Pi_{r,s}(x,y)$ with $r<n$ or $s<m$
such that $F(x,y)$ belongs to $\Pi_{r,s}(x,y)$.

\medskip

Let us note that these theorems refer to the degree of polynomials
in the power basis, so since now we will be using the Bernstein
basis some care will be needed. For the sake of clarity we will
illustrate all our results with a small example. Let

$$
\{\beta_0^{(4)}(t),\beta_1^{(4)}(t),\beta_2^{(4)}(t),\beta_3^{(4)}(t),\beta_4^{(4)}(t)\}
$$
be the (univariate) Bernstein basis of the space of polynomials of
degree less than or equal to $4$, where the Bernstein polynomials
are defined as follows,
$$
\beta_i^{(n)}(t) = {n \choose i} (1 - t)^{n-i} t^i, \qquad i = 0,
\ldots, n.
$$

and let us consider the algebraic curve given by the parametric
equation
$$
x(t)= {4 \beta_0^{(4)}(t)+ 4 \beta_1^{(4)}(t)+ 3 \beta_2^{(4)}(t)+
3 \beta_3^{(4)}(t)+ 7 \beta_4^{(4)}(t)\over  \beta_0^{(4)}(t)+
\beta_1^{(4)}(t)+  \beta_2^{(4)}(t)+  \beta_3^{(4)}(t)+ 3
\beta_4^{(4)}(t)}
$$
$$
y(t)= 2 \beta_0^{(4)}(t)+ 3 \beta_1^{(4)}(t)+ 3 \beta_2^{(4)}(t)+
3 \beta_3^{(4)}(t)+ 4 \beta_4^{(4)}(t).
$$

If we call $p(t) = u_1(t) - xv_1(t)$ and $q(t) = u_2(t) -
yv_2(t)$, their coefficients in the Bernstein basis are given by
$$
p_0 = 4 - x, p_1 = 4 - x, p_2 = 3 - x, p_3 = 3 - x, p_4 = 7 - 3x,
$$
and
$$
q_0 = 2 - y, q_1 = 3 - y, q_2 = 3 - y, q_3 = 3 - y, q_4 = 4 - y,
$$

However, let us observe that
$$
p(t) = 4 - x - 6t^2 + 8t^3 + (-2x + 1)t^4
$$
(a polynomial of degree $4$ in $t$), while

$$
q(t)= 2 - y + 4t - 6t^2 + 4t^3,
$$
a polynomial of degree $3$ in $t$.

Therefore, the polynomial defining the implicit equation will be a
polynomial belonging to the space $\Pi_{n,m}(x,y)$ with $n = 3$
and $m = 4$. We will use for that space the {\sl tensor-product
bivariate Bernstein basis} given by
$$
\{B_{ij}^{(n,m)}\} = \{\beta_i^{(n)}(x) \beta_j^{(m)}(y), i = 0,
\ldots, n ; j = 0, \ldots, m\}.
$$
\medskip

Finally we will recall, following [2], the algorithm for
constructing Bernstein-B\'ezout matrix of $p(t)$ and $q(t)$.
Although in [2] the coefficients of the polynomials are always
numbers, in our application we will construct the symbolic (i.e.
with the entries being polynomials in $x,y$) Bernstein-B\'ezout
matrix of $p(t)$ and $q(t)$ which we denote by $BS$. For the
reader's convenience, we present the algorithm written in {\sl
Maple} language:

\begin{verbatim}

for i from 1 to n do
   BS[i,1]:=(n/i)*(p[i]*q[0]-p[0]*q[i]);
od;

for j from 1 to n-1 do
   BS[n,j+1]:= (n/(n-j))*(p[n]*q[j]-p[j]*q[n])
od;

for j from 1 to n-1 do
   for i from 1 to n-1 do
     BS[i,j+1]:=(n^2/(i*(n-j)))*(p[i]*q[j]-p[j]*q[i])
     +((j*(n-i))/(i*(n-j)))*BS[i+1,j];
od;
od;

\end{verbatim}

Let us observe that if $m=n$, the resultant is the determinant of
the Bernstein-B\'ezout matrix, while -as a consequence of the
corresponding result for the B\'ezout resultant [16]- if $m>n$,
that determinant is equal to the resultant multiplied by the
factor $(\widetilde p_m(t))^{m-n}$, where $\widetilde p_m(t)$ is
the leading coefficient of $p(t)$ in the power basis.

So in our example, the determinant of $BS$ will be the implicit
equation we are looking for multiplied by the factor $(-2x + 1)$,
since the degree of $p$ is $4$ and the degree of $q$ is $3$ and
the coefficient of $t^4$ in $p$ is $(-2x + 1)$. In the following
section we will show how to compute the coefficients in the
bivariate tensor-product Bernstein basis of the implicit equation
(which will be a scalar multiple of the resultant computed by
using the approach of [12], where the equation is obtained in the
usual power basis).

\section{The interpolation process}

Since the expansion of the symbolic determinant is very time and
space consuming, our aim is to compute the polynomial defining the
implicit equation by means of Lagrange bivariate interpolation,
but using the Bernstein basis instead of the power basis. A good
introduction to the theory of interpolation can be seen in [5].

If we consider the {\it interpolation nodes} $(x_i,y_j)$ $(i = 0,
\cdots, n ; j = 0, \cdots, m)$ and the {\it interpolation space}
$\Pi_{n,m}(x, y)$, the interpolation problem is stated as follows:

Given $(n+1)(m+1)$ values
$$
f_{ij} \in K \qquad (i = 0, \cdots, n ; j = 0, \cdots, m)
$$
(the {\it interpolation data}), find a polynomial
$$
F(x,y) = \sum_{(i,j) \in I}c_{ij} \beta_i^{(n)}(x)
\beta_j^{(m)}(y)\in \Pi_{n,m}(x,y)
$$
(where $I$ is the index set $I = \{(i,j)|i = 0, \cdots, n ; j = 0,
\cdots, m\})$ such that
$$
F(x_i,y_j) = f_{ij} \qquad \forall \quad (i,j) \in I.
$$
If we consider for the interpolation space $\Pi_{n,m}(x, y)$ the
basis
$$
\{B_{ij}^{(n,m)}, ~~i = 0, \ldots, n ; j = 0, \ldots, m \} =$$ $$
\{ \beta_i^{(n)}(x) \beta_j^{(m)}(y), ~~ i = 0, \ldots, n ; j = 0,
\ldots, m \} = $$ $$ \{B_{00}^{(n,m)}, B_{01}^{(n,m)}, \cdots,
B_{0m}^{(n,m)}, B_{10}^{(n,m)}, B_{11}^{(n,m)}, \cdots,
B_{1m}^{(n,m)}, \cdots, $$ $$B_{n0}^{(n,m)}, B_{n1}^{(n,m)},
\cdots, B_{nm}^{(n,m)}\}
$$
with that precise ordering, and the interpolation nodes
$$
\{(x_i,y_j) | i = 0, \cdots, n; j = 0, \cdots, m\} =
$$
$$
 \{(x_0,y_0),(x_0,y_1), \cdots, (x_0,y_m),
$$
$$
(x_1,y_0),(x_1,y_1), \cdots,(x_1,y_m),\cdots,(x_n,y_0), \cdots,
(x_n,y_m)\},
$$
then the $(n+1)(m+1)$ interpolation conditions $F(x_i,y_j) =
f_{ij}$ can be written as a linear system

$$
A c = f,
$$
where the coefficient matrix $A$ is given by a Kronecker product
$$
B_x \otimes B_y,
$$
with
$$
B_x=((\beta_j^{(n)}(x_i)), ~~~i=0, \ldots,n; ~j=0,\ldots,n,
$$
$$
B_y=((\beta_j^{(m)}(y_i)), ~~~i=0, \ldots,m; ~j=0,\ldots,m,
$$
$$c=(c_{00}, \cdots, c_{0m}, c_{10}, \cdots, c_{1m}, \cdots, c_{n0}, \cdots, c_{nm})^T,$$
and
$$f=(f_{00}, \cdots, f_{0m}, f_{10}, \cdots, f_{1m}, \cdots, f_{n0}, \cdots, f_{nm})^T.$$
The Kronecker product $D \otimes E$ is defined by blocks as
$(d_{kl}E)$, with $D=(d_{kl})$.

For reasons which will be explained in Section $4$ we will select
as interpolation nodes $(x_i,y_j) =
({{i+1}\over{n+2}},{{j+1}\over{m+2}}) ~~(i = 0, \cdots, n; j = 0,
\cdots, m)$. In the general case we must avoid the value of $x_i$
for which the leading coefficient of $p(t)$ in the power basis
evaluates to $0$, and the value $y_j$ for which the leading
coefficient of $q(t)$ in the power basis evaluates to $0$.

In our example we have $n = 3$ and $m = 4$, and consequently $B_x$
will be the matrix
$$
B_x=\begin{pmatrix}
\frac{64}{125} & \frac{48}{125} & \frac{12}{125} & \frac{1}{125} \\
\frac{27}{125} & \frac{54}{125} & \frac{36}{125} & \frac{8}{125} \\
\frac{8}{125} & \frac{36}{125} & \frac{54}{125} & \frac{27}{125} \\
\frac{1}{125} & \frac{12}{125} & \frac{48}{125} & \frac{64}{125}
\end{pmatrix},
$$
and $B_y$ will be the matrix
$$
B_y=\begin{pmatrix}
\frac{625}{1296} & \frac{125}{324} & \frac{25}{216} & \frac{5}{324} & \frac{1}{1296} \\
\frac{16}{81} & \frac{32}{81} & \frac{8}{27} & \frac{8}{81} & \frac{1}{81}\\
\frac{1}{16} & \frac{1}{4} & \frac{3}{8} & \frac{1}{4} & \frac{1}{16} \\
\frac{1}{81} & \frac{8}{81} & \frac{8}{27} &
\frac{32}{81}&\frac{16}{81} \\
\frac{1}{1296} & \frac{5}{324} & \frac{25}{216} & \frac{125}{324}
& \frac{625}{1296}
\end{pmatrix}.
$$

As it is well known, since $B_x$ and $B_y$ are nonsingular
matrices the Kronecker product $B_x \otimes B_y$ will also be
nonsingular.

\bigskip

As for the generation of the interpolation data, let us remark
that they can be obtained without constructing the symbolic
Bernstein-B\'ezout matrix $BS$. That is to say, we can obtain each
interpolation datum by means of the evaluation of $p(t)$ and
$q(t)$ followed by the computation of the determinant of the
corresponding {\sl numerical} Bernstein-B\'ezout matrix $B$ making
use of the Bini-Gemignani algorithm which constructs (in $O(n^2)$
arithmetic operations) the Bernstein-B\'ezout matrix for the
evaluated polynomials.

In addition,  we must divide the value of the determinant by
$-2x_i+1$.

\medskip

An algorithm for solving linear systems with a Kronecker product
coefficient matrix is derived in a self-contained way (in a more
general setting) in [14]. For the case of the power basis
considered in [12], taking into account that every linear system
to be solved was a Vandermonde linear system, it was convenient to
use the Bj{\"o}rck-Pereyra algorithm [3, 9] to solve those linear
systems. For the Bernstein basis being used here, an appropriate
algorithm which takes advantage of the special properties of the
coefficient matrices $B_x$ and $B_y$ will be presented in Section
$4$.

In the general case, we must solve $n+1$ linear systems with the
same matrix $B_y$ and $m+1$ linear systems with the same matrix
$B_x$.

\section{Total positivity of $B_x$ and $B_y$}

From [4] we know that the Bernstein basis of the space of
polynomials of degree less than or equal to $n$ is a {\sl strictly
totally positive basis} on the open interval $(0,1)$, which
implies that all the {\sl collocation matrices}

$$
M = (\beta_j^{(n)}(t_i)), i,j = 0, \ldots, n
$$
with $t_0 < t_1 < \ldots < t_n$ in $(0,1)$ are {\sl strictly
totally positive}, i.e. all their minors are strictly positive. In
particular, due to our choice of the interpolation nodes the
matrices $B_x$ and $B_y$ are strictly totally positive matrices.

Making use of the results of [7, 8], we know that performing the
{\sl complete Neville elimination} on a strictly totally positive
matrix $A$ a {\sl bidiagonal factorization} of its inverse
$A^{-1}$ can be obtained, that is to say, we have

$$
A^{-1} = G_1 G_2 \ldots G_{n-1} D^{-1} F_{n-1} F_{n-2} \ldots F_1,
$$
where $D^{-1}$ is a diagonal matrix and $F_i$ and $G_i$ are
bidiagonal matrices.

So, after having obtained that factorization (with a computational
cost of $O(n^3)$ arithmetic operations), all the systems $A z = b$
with coefficient matrix $A$ can be solved (with a cost of $O(n^2)$
arithmetic operations) by performing the product

$$
G_1 G_2 \ldots G_{n-1} D^{-1} F_{n-1} F_{n-2} \ldots F_1 b.
$$

An early application of these ideas to solve structured linear
systems can be seen in [15], and a recent extension has been
presented in [6].

A detailed error analysis of Neville elimination, which shows the
advantages of this type of elimination for the class of totally
positive matrices, has been carried out in [1], and related work
for the case of Vandermonde linear systems can be seen in Chapter
$22$ of [10].

In our situation we must notice that the bidiagonal factorization
can be done in exact arithmetic, and the results of the
factorization can then be rounded if the subsequent computations
must be carried out in finite precision arithmetic.

After having obtained the bidiagonal factorization of the inverse
of $B_y$, the solution of the linear system $B_y z = b$ can be
obtained in $O(n^2)$ arithmetic operations by computing the
product

$$
G_1 G_2 \ldots G_{n-1} D^{-1} F_{n-1} F_{n-2} \ldots F_1 b,
$$
and analogously for the linear systems with coefficient matrix
$B_x$ [6].

In our example, the coefficients of the desired implicit equation
in the tensor-product bivariate Bernstein basis (using the
lexicographical ordering we are considering) are:

$$
(25264/27, 66256/81, 167852/243, 45652/81, 36137/81, 15728/27,
125312/243,
$$
$$
320120/729, 29164/81, 69421/243, 29440/81, 79024/243, 203228/729,
$$
$$
18580/81, 14761/81, 2048/9, 16640/81, 14336/81, 3940/27, 9391/81).
$$

\section{Computational complexity}

In this section we will briefly examine the computational
complexity of our algorithm in terms of arithmetic operations. In
view of the algorithm, we must solve $n+1$ systems of order $m+1$
with the same matrix $B_y$ and $m+1$ systems of order $n+1$ with
the same matrix $B_x$.

\medskip

The factorization of the inverse of a matrix of order $n$ by means
of complete Neville elimination takes $O(n^3)$ operations, but
that factorization is used for solving all the systems with the
same matrix, so each of the remaining systems can be solved with
$O(n^2)$ operations.

For the sake of clarity in the comparison, we will consider here
the case $m=n$. Then, the interpolation part of the algorithm has
computational complexity $O(n^3)$. Let us observe that in this
situation, if we solve the linear system $Ac=f$ of order $(n+1)^2$
by means of Gaussian elimination, without taking into account the
special structure of the matrix, we have computational complexity
$O(n^6)$. Moreover, using the approach we are describing, there is
no need of constructing the matrix $A$, which implies an
additional saving in computational cost.

Let us remark that, since the construction of the numerical
Bernstein-B\'ezout matrix requires $O(n^2)$ arithmetic operations
and the complexity of the computation of each determinant is
$O(n^3)$, the generation of the interpolation data has a
computational complexity of $O(n^5)$. Therefore with our approach,
which exploits the Kronecker product structure, the whole process
has complexity $O(n^5)$, while using Gaussian elimination it would
be $O(n^6)$.

It is worth noting that the main cost of the process corresponds
to the generation of the interpolation data, and not to the
computation of the coefficients of the interpolating polynomial.
So, the main effort to reduce the computational cost must be
focused on that stage. In this sense, an interesting issue would
be to take advantage of the {\sl displacement structure} of the
Bernstein-B\'ezout matrices [2, 11] to develop an algorithm with
complexity $O(n^2)$ for computing each determinant.

\medskip

{\sl Remark}. Finally, let us observe that all the linear systems
with matrix $B_y$ can be solved simultaneously, and the same can
be said of the systems with matrix $B_x$. Therefore the algorithm
exhibits a high degree of {\it intrinsic parallelism}. This
parallelism is also present in the computation of the
interpolation data since we can compute simultaneously the
determinants involved in this process.

\bigskip

{\bf Acknowledgements}

This research has been partially supported by Spanish Research
Grant BFM 2003-03510 from the Spanish Ministerio de Ciencia y
Tecnolog{\'\i}a.

\end{document}